\documentclass[leqno,12pt]{article}
\usepackage{latexsym,amsmath,amssymb,amscd}

\title{\uppercase{On the Tolman-Oppenheimer-Volkoff-de Sitter equation} }
\author{Tetu Makino}
\date{\today}
\begin{document}
\maketitle
\footnote{ 2010 {\it Mathematical Subject Classification.} Primary 34A34;
Secondary 37N20, 76N10, 83C05, 83C20.
}
\footnote{{\it Key words and phrases.} Einstein equations, de Sitter models, Cosmological constant,
Spherically symmetric solutions, Vacuum boundary 
}

\begin{abstract}
Spherically symmetric static solutions of the Einstein
equations with a positive cosmological
constant for the energy-momentum tensor of a barotropic
perfect fluid are governed by the Tolman-Oppenheimer-Volkoff-de Sitter equation. Sufficient conditions for existence of solutions with finite radii are given. The interior metric of the solution is connected with the Schwarzschild-de Sitter metric on the exterior vacuum region. The analytic property of the solutions at the vacuum boundary is investigated. 
\end{abstract}

\newtheorem{Lemma}{Lemma}
\newtheorem{Proposition}{Proposition}
\newtheorem{Theorem}{Theorem}
\newtheorem{Definition}{Definition}

\section{Introduction}

We consider a static and spherically symmetric metric
$$
ds^2=g_{\mu\nu}dx^{\mu}dx^{\nu}=e^{2F(r)}c^2dt^2-e^{2H(r)}dr^2-
r^2(d\theta^2+\sin^2\theta d\phi^2)
$$
which satisfies the Einstein-de Sitter equations
$$
R_{\mu\nu}-\frac{1}{2}g_{\mu\nu}R-
\Lambda g_{\mu\nu}=\frac{8\pi G}{c^4}T_{\mu\nu}, 
$$
for the energy-momentum tensor of a perfect fluid
$$
T^{\mu\nu}=(c^2\rho +P)U^{\mu}U^{\mu}-Pg^{\mu\nu}.
$$
Here
$R_{\mu\nu}$ is the Ricci tensor associated with the metric
$g_{\mu\nu}dx^{\mu}dx^{\nu}$, $R=g^{\alpha\beta}R_{\alpha\beta}$ is the 
scalar curvature, and  $c, G$ are positive constants, the speed of light, the gravitational
 constant. $\Lambda$ is the cosmological constant which is supposed 
to be positive. $\rho$ is the mass density, $P$ is the pressure and $U^{\mu}$ is the four-velocity. See \cite[\S 111]{LandauL}.

Historically speaking, the cosmological constant $\Lambda$ of
the above Einstein-de Sitter equations was introduced by A. Einstein
\cite{Einstein}, 1917,
and was discussed soon by W. de Sitter \cite{deSitter}. Although it was introduced for a static universe, it was not necessary for an expanding universe.
So, later Einstein wrote to H. Weyl ``If there is no quasi-static world,
then away with the cosmological term! " on May 23, 1923, and
finally he rejected the cosmological term in \cite{Einstein1931}, 1931, for the reason that it is not necessary to explain 
the Hubble's report on the redshifts of galaxies showing the expansion of the universe. (See
\cite[\S 15e]{Pais}.) However, although the original motivation for 
introducing the cosmological term disappeared, its status in the 
cosmological theories remained and it revived with new meanings in the recent development of the theories and observations. For the details, see the review \cite{Carroll}\\ 

For the static spherically symmetric metric, the Einstein-de Sitter equations are reduced to
\begin{equation}
\frac{dm}{dr}=4\pi r^2\rho, \quad 
\frac{dP}{dr}=-(\rho+P/c^2)
\frac{\displaystyle G\Big(m+\frac{4\pi r^3}{c^2}P\Big)-\frac{c^2\Lambda}{3}r^3}{r^2\Big(\displaystyle1-\frac{2Gm}{c^2r}-\frac{\Lambda}{3}r^2\Big)}. \
\label{TOV1}
\end{equation}
The coefficients of the metric are given by $$e^{2F(r)}=\kappa_+e^{-2u(r)/c^2},
$$
$u(r), \kappa_+$ being the function and the constant specified later, and
$$e^{-2H(r)}=1-\frac{2Gm(r)}{c^2r}-\frac{\Lambda}{3}r^2.
$$

When $\Lambda=0$, the equations (\ref{TOV1}) turns out to be
\begin{equation}
\frac{dm}{dr}=4\pi r^2\rho, \quad \frac{dP}{dr}=-(\rho+P/c^2)
\frac{\displaystyle G\Big(m+\frac{4\pi r^3}{c^2}P\Big)}{r^2\Big(\displaystyle1-\frac{2Gm}{c^2r}\Big)}, \label{TOV8}
\end{equation}
and this (\ref{TOV8}) is called the Tolman-Oppenheimer-Volkoff equation. It was derived in \cite{OppenheimerV}, 1939.
Therefore we shall call (\ref{TOV1}) with $\Lambda>0$ {\bf the
Tolman-Oppenheimer-Volkoff-de Sitter equation}. In this article we investigate this Tolman-Oppenheimer-Volkoff-de Sitter equation 
(\ref{TOV1}).  Throughout this article we keep the following\\

{\bf Assumption}  {\it The pressure $P$ is
 a given function of the density $\rho>0$
such that $0<P$ and $0<dP/d\rho <c^2$ for
$\rho>0$ and $P \rightarrow 0$ as $\rho \rightarrow +0$. 
Moreover we assume that there are positive constants $A, \gamma$ and an analytic
function $\Omega$ on a neighborhood of $[0, +\infty)$ such that
$\Omega(0)=1$ and
$$P=A\rho^{\gamma}\Omega(A\rho^{\gamma-1}/c^2).
$$
We assume that $1<\gamma<2$.}\\

Actually we are keeping in mind the equation of state for neutron
stars:
$$P=Kc^5\int_0^{\zeta}
\frac{q^4dq}{\sqrt{1+q^2}},\quad
\rho=3Kc^3\int_0^{\zeta}\sqrt{1+q^2}q^2dq, $$
$K$ being a positive constant. 
 See \cite{2015KJM}, \cite[p. 188]{ZeldovichN}.\\

The Tolman-Oppenheimer-Volkoff-de Sitter equation has already
been systematically investigated from the physical point of view
by C. G. B\"{o}hmer, \cite{Boehmer1}, \cite{Boehmer2}. However
in the study by C. G. B\"{o}hmer there is supposed to exist
a positive density $\rho_b$, so called the `boundary density',
at which the pressure $P$ vanishes, that is,
$P>0 \Leftrightarrow \rho>\rho_b$, and cosmological constants
$\Lambda$ satisfying $\displaystyle \Lambda <\frac{4\pi G}{c^2}\rho_b$
are considered. This situation is not treated in this article.\\

For the sake of notational conventions, we shall denote
\begin{equation}
\kappa(r,m):=1-\frac{2Gm}{c^2r}-\frac{\Lambda}{3}r^2, \label{kappa}
\end{equation}
and
\begin{equation}
Q(r,m,P):=G\Big(m+\frac{4\pi r^3}{c^2}P\Big)-\frac{c^2\Lambda}{3}r^3, \label{Q}
\end{equation}
so that the second equation of (\ref{TOV1}) reads
$$\frac{dP}{dr}=-(\rho +P/c^2)\frac{Q(r, m, P)}{r^2\kappa(r,m)}.$$\\

We consider the equation (\ref{TOV1}) on the domain
$$\mathcal{D}=\{(r,m,P) | 0<r, |m|<+\infty,
0<\rho, 0<\kappa(r,m) \}.
$$\\

Now we have a solution germ $(m(r), P(r))$
at $r=+0$, given the central density $\rho_c>0$ with $P_c:=P(\rho_c)$,
such that
\begin{align}
m&=\frac{4\pi}{3}\rho_cr^3+O(r^5), \label{TOV2} \\
P&=P_c-
(\rho_c+P_c/c^2)
\Big(4\pi G(\rho_c+3P_c/c^2)-c^2\Lambda\Big)
\frac{r^2}{6}+O(r^4) \nonumber
\end{align}
as $r\rightarrow +0$. 
Proof is the same as that of \cite[Proposition 1]{2015KJM}. See
\cite[\S 2, pp. 57-58]{1998JMKU}.

We are interested in the prolongation of the solution germ
to the right as long as possible in the domain $\mathcal{D}$.
Actually the prolongation is unique, since the right-hand sides
of (\ref{TOV1}) are analytic functions of
$r, m, P$ in $\mathcal{D}$. See \cite[Chap.1, Sec.5]{CoddingtonL}.

Especially we want to have a sufficient condition for that
the prolongation turns out to be `monotone-short' in the following
sense:

\begin{Definition}
A solution $(m(r), P(r)), 0<r<r_+,$ of
(\ref{TOV1}) is said to be {\bf monotone-short} if
$r_+<\infty$, $dP/dr <0$ for $0<r<r_+$ and
$P\rightarrow 0$ as $r\rightarrow r_+$
and if $\kappa_+>0$ and $Q_+>0$, where
\begin{align*}
\kappa_+&:=\lim_{r\rightarrow r_+-0}\kappa(r, m(r))=1-\frac{2Gm_+}{c^2r_+}-\frac{\Lambda}{3}r_+^2, \\
Q_+&:=\lim_{r\rightarrow r_+-0} Q(r, m(r), P(r))=Gm_+-\frac{c^2\Lambda}{3}r_+^3,
\end{align*}
with $$m_+:=\lim_{r\rightarrow r_+-0} m(r).$$
\end{Definition}

Since the solution germ behaves as (\ref{TOV2}) as $r\rightarrow +0$,
we assume 
\begin{equation}
\Lambda<\frac{4\pi G}{c^2}(\rho_c+3P_c/c^2) \label{TOV3}
\end{equation}
in order that $dP/dr< 0$ at least for $0<r\ll 1$. 

{\bf Remark.}  If  the condition (\ref{TOV3}) does not hold
but the equality
\begin{equation}
\Lambda = \frac{4\pi G}{c^2}(\rho_c+3P_c/c^2) \label{L}
\end{equation}
holds exactly, then the solution of (\ref{TOV1}), (\ref{TOV2}) 
turns out to be $(m(r), P(r))=(4\pi \rho_c r^3/3, \  P_c),
0<r<\sqrt{3/L}$, where
$\displaystyle L:=\frac{8\pi G}{c^2}\rho_c+\Lambda.
$
In this very special case, the metric is reduced to
$$
ds^2=c^2dt^2-\Big(1-\frac{L}{3}r^2\Big)^{-1}dr^2-r^2(d\theta^2+\sin^2\theta d\phi^2)
$$
after a suitable change of the scale of $t$.
The space at $t=\mbox{Const.}$ is isometric to the half of the compact
3-dimensional hypersphere with radius $\sqrt{3/L}$
embedded in the 4-dimensional Euclidean space
$\mathbb{R}^4=\{ (\xi^1,\xi^2,\xi^3,\xi^4)\}$
through
$$\xi^1=r\sin\theta\cos\phi, \quad \xi^2=r\sin\theta\sin\phi, \quad \xi^3=r\cos\theta,
\quad \xi^4=\displaystyle\sqrt{\frac{3}{L}-r^2}.$$
Thus the horizon $r=\sqrt{3/L}-0$ is merely
apparent. Of course the condition (\ref{L}) is highly unstable.
This is nothing but the `Einstein's steady state
inverse (1917)' proposed in \cite{Einstein}. 
On the other hand, suppose that the inequality
\begin{equation}
\Lambda >\frac{4\pi G}{c^2}(\rho_c+3P_c/c^2) \label{LL}
\end{equation}
holds. Then the solution germ
satisfies $dP/dr >0$ for $ 0<r \ll 1$. However it is possible that
$dP/dr$ become negative when prolonged to the right. 
This fact will be shown later. $\square$\\

We introduce the variable $u$ by
$$u:=
\int_0^{\rho}
\frac{dP}{\rho+P/c^2}.
$$
Then we see
\begin{align*}
u&=\frac{\gamma A}{\gamma-1}\rho^{\gamma-1}\Omega_u(A\rho^{\gamma-1}/c^2),\\
\rho&=A_1u^{\frac{1}{\gamma-1}}\Omega_{\rho}(u/c^2), \\
P&=AA_1^{\gamma}u^{\frac{\gamma}{\gamma-1}}\Omega_P(u/c^2),
\end{align*}
where $\Omega_u, \Omega_{\rho}, \Omega_P$ are 
analytic functions on a neighborhood of $[0,+\infty[$
such that $\Omega_u(0)=\Omega_{\rho}(0)=\Omega_P(0)=1$
and$\displaystyle A_1:=\Big(\frac{\gamma-1}{\gamma A}\Big)^{\frac{1}{\gamma-1}}$. The functions $\Omega_u, \Omega_{\rho}, \Omega_P$ are 
depending upon only $\gamma$ and the function $\Omega$.
In fact we take
\begin{align*}
\Omega_u(\zeta)&=
\frac{1}{\zeta}\int_0^{\zeta}
\frac{\Omega(\zeta')+\frac{\gamma-1}{\gamma}\zeta'D\Omega(\zeta')}{1+\zeta'\Omega(\zeta')}d\zeta', \\
\zeta&=\frac{\gamma-1}{\gamma}\eta \Omega_{\rho}(\eta) \Leftrightarrow 
\eta =\frac{\gamma}{\gamma-1}\zeta\Omega_u(\zeta), \\
\Omega_P(\eta)&=
\Omega(\zeta)\Omega_u(\zeta)^{-\frac{\gamma}{\gamma-1}} \quad\mbox{with}\quad
\zeta=\frac{\gamma-1}{\gamma}\eta\Omega_{\rho}(\eta)
\end{align*}
Let us fix a small positive number $\delta_{\Omega}$ such that
these functions are defined and analytic on a neighborhood of 
$[-\delta_{\Omega}, +\infty[$. 
We put
$$
u_c:=\int_0^{\rho_c}\frac{dP}{\rho+P/c^2}=\frac{\gamma A}{\gamma-1}\rho_c^{\gamma-1}
\Omega_u(A\rho_c^{\gamma-1}/c^2).
$$

\section{Main result}

We claim

\begin{Theorem}
Suppose that $6/5<\gamma<2$. Then there exists a positive number
$\epsilon_0(\leq 1)$ depending upon only $\gamma$ and the function
$\Omega$ such that if
\begin{equation}
u_c\leq c^2\epsilon_0,\qquad
\Lambda\leq \frac{4\pi}{c^2}G\Big(\frac{\gamma-1}{\gamma A}\Big)^{\frac{1}{\gamma-1}}
(u_c)^{\frac{1}{\gamma-1}}\epsilon_0,\label{TOV4}
\end{equation}
then the prolongation of the solution germ with (\ref{TOV2})
to the right turns out to be monotone-short.
\end{Theorem}

This can be considered as the de Sitter version of the result
by A. D. Rendall and G. B. Schmidt, \cite{RendallS}.\\

Let us sketch the proof. 
Using the variable $u$, we can write the Tolman-Oppenheimer-Volkoff-de Sitter equation (\ref{TOV1}) as
\begin{align}
\frac{dm}{dr}&=
4\pi r^2A_1(u_{\sharp})^{\frac{1}{\gamma-1}}\Omega_{\rho}(u/c^2), \label{TOV5} \\
\frac{du}{dr}&=
-\frac{G\Big(\displaystyle m+\frac{4\pi}{c^2}r^3AA_1^{\gamma}(u_{\sharp})^{\frac{\gamma}{\gamma-1}}\Omega_P(u/c^2)\Big)-\frac{c^2\Lambda}{3}r^3}{r^2\Big(\displaystyle 1-\frac{2Gm}{c^2r}-\frac{\Lambda}{3}r^2\Big)}. \nonumber
\end{align}
Here $(u_{\sharp})$ stands for $\max \{u, 0\}$. Since we are assuming that $1<\gamma<2$, which
implies $\displaystyle \mu:=\frac{1}{\gamma-1}>1, \frac{\gamma}{\gamma-1}=\mu+1>2$,
we can see that the functions $u\mapsto (u_{\sharp})^{\mu}, (u_{\sharp})^{\mu+1}$ are of class
$C^1(\mathbb{R})$. Keeping in mind it, we consider that the domain of
the equation(\ref{TOV5}) is
$$
\mathcal{D}_u=\{(r,m,u)| 0<r, |m|<\infty,
-\delta_{\Omega}<u/c^2<+\infty,
\kappa >0\}.
$$

Let us perform the homologous transformation of the variables
$$r=aR, \qquad m=a^3b^{\frac{1}{\gamma-1}}\cdot 4\pi A_1 M,
\qquad u=bU,
$$
where $a,b$ are positive parameters. We take $b=u_c$ and $a$ which satisfies
$$4\pi GA_1a^2b^{\frac{2-\gamma}{\gamma-1}}=1.$$ Let us write
$$\lambda :=\frac{c^2}{4\pi GA_1}\Lambda,
\qquad \alpha:=b/c^2=u_c/c^2,
\qquad \beta:=b^{-\frac{1}{\gamma-1}}\lambda=
\frac{c^2}{4\pi GA_1}(u_c)^{-\frac{1}{\gamma-1}}\Lambda.
$$
Then the system (\ref{TOV5}) turns out to
\begin{align}
\frac{dM}{dR}&=R^2(U_{\sharp})^{\frac{1}{\gamma-1}}
\Omega_{\rho}(\alpha U), \label{TOV6} \\
\frac{dU}{dR}&=-\frac{1}{R^2}
\frac{\displaystyle\Big(M+\frac{\gamma-1}{\gamma}\alpha R^3(U_{\sharp})^{\frac{\gamma}{\gamma-1}}\Omega_P(\alpha U)
-\frac{1}{3}\beta  R^3\Big)}{\displaystyle\Big(1-2\alpha\frac{M}{R}-
\frac{1}{3}\alpha\beta R^2\Big)}. \nonumber
\end{align}
Here $(U_{\sharp})$ stands for $\max\{U, 0\}$.
The domain of the system (\ref{TOV6}) should be
$$\mathcal{D}_U=\{(R, M, U) |
0<R, |M|<\infty,
-\delta_{\Omega}<U<2, \kappa>0\},
$$
where, of course,
$$\kappa=1-2\alpha\frac{M}{R}-
\frac{1}{3}\alpha\beta R^2.$$
Let us concentrate ourselves to $0\leq \alpha\leq 1, 0\leq \beta\leq 1$.
We are considering a solution germ $(M(R), U(R))$ at $R=+0$
which satisfies
\begin{align}
M(R)&=\Omega_{\rho}(\alpha)\frac{R^3}{3}+O(R^5), \label{MU0} \\
U(R)&=1-
\Big(\Omega_{\rho}(\alpha)+\frac{3\gamma}{\gamma-1}
\alpha\Omega_P(\alpha)-\beta\Big)\frac{R^2}{6}+O(R^4). \nonumber
\end{align}

We claim
\begin{Proposition}
There is a positive number $R_0$ which depends upon only
$\gamma$ and the function $\Omega$ such that
$(M(R), U(R))$ exists and satisfies $dU/dR<0$ on 
$0<R \leq R_0$ and $(M(R_0), U(R_0))$ depends continuously
on $\alpha, \beta \in [0,1]$.
\end{Proposition}
Proof is standard, and done by converting the system of differential
equations 
(\ref{TOV6}) to a system of integral equations
under the condition (\ref{MU0}) as
\begin{align}
q(R)&=\frac{3}{R^3}\int_0^R
U(R')^{\frac{1}{\gamma-1}}
\frac{\Omega_{\rho}(\alpha U(R'))}{\Omega_{\rho}(\alpha)}
R'^2dR', \label{MUinte} \\
U(R)&=1-\int_0^R
\frac{\frac{1}{3}\Omega_{\rho}(\alpha)q(R')+\frac{\gamma-1}{\gamma}\alpha
U(R')^{\frac{\gamma}{\gamma-1}}
\Omega_P(\alpha U(R'))-\frac{1}{3}\beta}{1-2\alpha\frac{1}{3}\Omega_{\rho}(\alpha)q(R')R'^2-\frac{1}{3}\alpha\beta R'^2}R'dR'.
\nonumber
\end{align}
Then, taking $\delta$ sufficiently small uniformly on $\alpha, \beta$, we see that the mapping $(q, U) \mapsto (\tilde{q}, \tilde{U})$, which is
the right-hand side of the (\ref{MUinte}), is a contraction from
$\mathfrak{F}=\{(q,U) \in C[0,\delta] | 0\leq q \leq C_q, \frac{1}{2}\leq U\leq 2\}$ into 
itself with respect to a suitable functional distance,
where $$C_q:=\max\Big\{ U^{\frac{1}{\gamma-1}}\frac{\Omega_{\rho}(\alpha U)}{\Omega_{\rho}(\alpha)}\ \Big|\  
\frac{1}{2} \leq U\leq 2, \quad 0\leq \alpha \leq 1\Big\}.$$ For the details see \cite[pp. 57-58]{1998JMKU}.\\

Let us come back to the proof of Theorem 1. The right-hand side of the system (\ref{TOV6}) depends continuously
on $\alpha, \beta$, and tends to 
$(R^2(U_{\sharp})^{\frac{1}{\gamma-1}}, -M/R^2)^T$ as 
$\alpha \rightarrow 0, \beta\rightarrow 0$.
The limit system
$$
\frac{dM}{dR}=R^2(U_{\sharp})^{\mu}, \qquad
\frac{dU}{dR}=-\frac{M}{R^2}
$$
is nothing but the Lane-Emden equation
\begin{equation}
-\frac{1}{R^2}\frac{d}{dR}\Big(R^2\frac{dU}{dR}\Big)=(U_{\sharp})^{\mu}.
\label{LE}
\end{equation}
Since we are assuming that $6/5<\gamma<2$, say,
$1<\mu<5$, the solution $U=\bar{U}(R)$
with $\bar{U}(0)=1$ of the Lane-Emden equation (\ref{LE}) is short, that is,
$0<\bar{U}(R), d\bar{U}/dR<0$ for
$0<R<\xi_1=\xi_1(\gamma)$ and $\bar{U}(\xi_1)=0$. See \cite{Chandrasekhar}, \cite{JosephL}. 
Of course we consider
$$
\bar{U}(R)=\Big(R^2\frac{d\bar{U}}{dR}\Big)_{R=\xi_1}
\Big(\frac{1}{\xi_1}-\frac{1}{R}\Big)
$$
harmonically on $R\geq\xi_1$.  

Thanks to Proposition 1, if $\epsilon_0$ is sufficiently small and
if
\begin{equation}
\alpha\leq\epsilon_0 \qquad \mbox{and} \qquad
\beta \leq \epsilon_0, \label{TOV7}
\end{equation}
then $U=U(R)$ exists and remains near to the orbit of
$U=\bar{U}(R)$ on $R_0\leq R\leq \xi_1+\delta_R$,
$\delta_R$ being small so that $$
-\frac{\delta_{\Omega}}{2}\leq \bar{U}(\xi_1+\delta_R)<0.$$
This is nothing but a direct application of
\cite[Theorem 7.4]{CoddingtonL}.
Note that $\max\{\Omega_P(\eta) | -\delta_{\Omega}\leq
\eta\leq 2\}$ depends upon only $\gamma$ and the function
$\Omega$, and we have $-\delta_{\Omega}\leq \eta=\alpha U\leq 2$
provided that $\alpha\leq\epsilon_0\leq 1$ and 
$-\delta_{\Omega}<U<2$. Especially if $U(\xi_1+\delta_R)<0$,
then the radius $R_+$ of $U(R)$ should be found 
in the interval $]0, \xi_1+\delta_R[$. This completes the proof of Theorem 1,
since the condition (\ref{TOV7}) is nothing but (\ref{TOV4}).\\

Note that (\ref{TOV3}) follows from (\ref{TOV4}) if $\epsilon_0$ is
sufficiently small, since
$$\rho+3P/c^2=\Big(\frac{\gamma-1}{\gamma A}\Big)^{\frac{1}{\gamma-1}}
u^{\frac{1}{\gamma-1}}
\Omega_{\rho+3P/c^2}(u/c^2),
$$
where
$$\Omega_{\rho+3P/c^2}(\eta):=\Omega_{\rho}(\eta)
+
3\frac{\gamma-1}{\gamma}
\eta\Omega_P(\eta)
$$
is a function depending upon only $\gamma$ and
the function $\Omega$ and $\Omega_{\rho+3P/c^2}(0)=1$
so that we can assume that
$\min\{\Omega_{\rho+3P/c^2}(\eta)|
-\delta_{\Omega}\leq\eta\leq 2\}>\epsilon_0$.\\

{\bf Remark.}  For the existence of $u_c$ satisfying (\ref{TOV4}) it is necessary that $\Lambda$ enjoys
$$\Lambda \leq 4\pi c^{\frac{2(2-\gamma)}{\gamma-1}}G
\Big(\frac{\gamma-1}{\gamma A}\Big)^{\frac{1}{\gamma-1}}\epsilon_0^{\frac{\gamma}{\gamma-1}}.
$$
Thus one may ask whether the real value of the cosmological constant of our universe
satisfies it or not. But this question is not theoretical but experimental-observational and numerical.
To answer to it is not a business of such a poor mathematician as the author
of this article. $\square$\\

Even if $\gamma \leq 6/5$, the solution with the central density
$\rho_c$ of the Tolman-Oppenheimer-Volkoff equation,
that is, (\ref{TOV1}) with $\Lambda=0$, or (\ref{TOV8}),
can be short, if $\rho_c$ is large and the function $P(\rho)$ is very much different
from the exact $\gamma$-law for large $\rho$. 
For a sufficient condition for solutions to be short, see
 \cite[Proposition 3]{2015KJM}, the proof of 
\cite[Theorem 1]{1998JMKU}.
Therefore we 
can consider such a case, supposing that $1<\gamma<2$ and
the solution $(m,P)=(m^0(r), P^0(r)), 0<r<r_+^0,$ of (\ref{TOV8}) with
the same central density $\rho_c$ satisfies $P^0(r)\rightarrow 0$
as $r\rightarrow r_+^0-0 $, with $r_+^0$ being finite. Then the associated 
$(m,u)=(m^0(r), u^0(r)), 0<r<r_+^0,$
satisfies (\ref{TOV5}) with $\Lambda=0$, that is,
\begin{align}
\frac{dm}{dr}&=
4\pi r^2A_1(u_{\sharp})^{\frac{1}{\gamma-1}}\Omega_{\rho}(u/c^2), \label{TOV9} \\
\frac{du}{dr}&=
-\frac{G\Big(\displaystyle m+\frac{4\pi}{c^2}r^3AA_1^{\gamma}(u_{\sharp})^{\frac{\gamma}{\gamma-1}}\Omega_P(u/c^2)\Big)}{r^2\Big(\displaystyle 1-\frac{2Gm}{c^2r}\Big)}. \nonumber
\end{align}
In order to extend $(m^0(r), u^0(r))$ onto
$r\geq r_+^0$, we put
\begin{align*}
m^0(r)&=m_+ ^0(:=m^0(r_+^0)), \\
u^0(r)&=\frac{c^2}{2}\Big(
\log\Big(1-\frac{2Gm_+^0}{c^2r_+^0}\Big)
-\log\Big(1-\frac{2Gm_+^0}{c^2r}\Big)\Big),
\end{align*}
for $r\geq r_+^0$. Then the extended $(m^0(r),u^0(r))$
satisfies (\ref{TOV9}) on $0<r<+\infty$.
Since the right-hand side of (\ref{TOV5}) tends to that of (\ref{TOV9})
as $\Lambda \rightarrow 0$,
the solution of (\ref{TOV5}) under consideration
exists and remains in a neighborhood of the 
orbit $(m^0(r),u^0(r))$ on 
$0<r\leq r_+^0+\delta_r$, $\delta_r$ being a sufficiently small positive
number, provided that $\Lambda$ is sufficiently small.
Thus we have
\begin{Theorem}
Suppose that the solution of the Tolman-Oppenheimer-Volkaff equation
(\ref{TOV8}) with central density $\rho_c$ is short. Then there exists 
a small positive number $\epsilon_1$ such that,
if $\Lambda \leq\epsilon_1$, the solution germ
of the Tolman-Oppenheimer-Volkoff-de Sitter equation
(\ref{TOV1}) with the central density $\rho_c$
has a monotone-short prolongation.
\end{Theorem}

We should note that $\epsilon_1$ may depend not only upon
$\gamma$ and the function $\Omega$ but also
upon $A, c, G$ and $\rho_c$. In contrast with Theorem 1, we have no hope to specify 
the manner of dependence. 

\section{Monotonicity}

Here let us give a remark on the monotonicity of the
solution of the Tolman-Oppenheimer-Volkoff-de Sitter
equation (\ref{TOV1}).\\

When we studied the Tolman-Oppenheimer-Volkoff equation,
that is, (\ref{TOV1}) with $\Lambda=0$, or (\ref{TOV8}), we see
that if $]0,r_+[, r_+\leq +\infty,$ is the right maximal interval of
existence of the solution in the domain $\mathcal{D}$, then
$dP/dr <0$ for $0<r<r_+$ and $P \rightarrow 0$ as
$r\rightarrow r_+-0$. Proof is given in \cite{2015KJM}, \cite{1998JMKU}.
In other words, we can say on the Tolman-Oppenheimer-Volkoff
equation that, if the prolongation of the solution germ
under consideration is short, it is necessarily monotone-short.
However it is not the case on the Tolman-Oppenheimer-Volkoff-de Sitter
equation with $\Lambda>0$. Even if $dP/dr<0$
for $0<r\ll 1$ under the assumption (\ref{TOV3}), 
$dP/dr$ may turn out to be positive during the prolongation.
Let us show it.

In order to fix the idea, we suppose $u_c=1$, and put
$$r=aR,\quad m=a^3\cdot 4\pi A_1M, \quad
4\pi GA_1a^2=1,
\quad u=U, 
\quad \lambda=\frac{c^2}{4\pi GA_1}\Lambda.
$$
Then the system (\ref{TOV1}) is reduced to
\begin{align}
\frac{dM}{dR}&=R^2U^{\mu}\Omega_{\rho}(U/c^2), \label{TOV11} \\
\frac{dU}{dR}&=-\frac{1}{R^2}
\Big(M+\frac{\gamma-1}{\gamma}\frac{R^3}{c^2}U^{\mu +1}
\Omega_P(U/c^2)-
\frac{\lambda}{3}R^3\Big)\times \nonumber \\
&\times
\Big(1-\frac{2M}{c^2R}-\frac{\lambda}{3c^2}R^2\Big)^{-1}
. \nonumber
\end{align}
Here $\mu:=1/(\gamma-1)$.

The right-hand side of the system (\ref{TOV11}) depends continuously upon the speed of light $c$ and tends to 
$\displaystyle\Big(R^2U^{\mu}, -\frac{1}{R^2}\Big(M-\frac{\lambda}{3}R^3\Big)\Big)^T$ as $c\rightarrow \infty$. This non-relativistic limit 
equation can be written as
\begin{equation}
-\frac{1}{R^2}\frac{d}{dR}\Big(R^2\frac{dU}{dR}\Big)=U^{\mu}-\lambda.
\label{TOV12}
\end{equation}
In this situation we are assuming that $\Lambda$ depends upon $c$
and
$c^2\Lambda/(4\pi GA_1)$ tends to $\lambda$.
Since (\ref{TOV12}) is the Lane-Emden equation when $\lambda=0$,
we shall call it `{\bf the Lane-Emden-de Sitter equation}' supposing that
$\lambda>0$.

Although we are supposing $1<\mu<+\infty (
\Leftrightarrow 1<\gamma<2)$, we observe the limiting case
$\mu=1 (\Leftrightarrow \gamma=2)$. Then the equation (\ref{TOV12})
is linear and the solution $U=\hat{U}(R)$ with
$\hat{U}(0)=1$ is given by 
$$
\hat{U}(R)=\lambda+(1-\lambda)\frac{\sin R}{R}
$$
explicitly. We have
$$\hat{U}(R)=1-\frac{1-\lambda}{6}R^2+O(R^4)
$$
as $R \rightarrow +0$, and the condition
(\ref{TOV3}) reads $\lambda<1$.
Suppose that $\frac{1}{2}\leq\lambda <1$. Then we find that
$$
\frac{d\hat{U}}{dR}=\frac{1-\lambda}{R}\Big(\cos R-\frac{\sin R}{R}\Big)
$$ turns out to be positive for $3\pi/2 <R<2\pi$ and so on, while
$\hat{U}(R) >0$ exists and oscillates on $0<R<+\infty$, and converges to
$\lambda$ as $R\rightarrow +\infty$. Therefore,  
this explicit example tells us that, if
$\gamma$ is near to 2, $c$ is sufficiently large, and
$c^2\Lambda/(4\pi GA_1)$ is near to a number $\lambda$ in
the interval $[\frac{1}{2},1[$, the behavior of the solution under consideration may be similar, that is, the prolongation of the solution germ 
with $u_c=1$ is not monotone. On the other hand, the condition
(\ref{LL}) reads $\lambda >1$. Then $d\hat{U}/dR >0$ for $0<R\ll 1$ but
$d\hat{U}/dR$ become negative and $\hat{U}(R)$ oscillates and tends to
the limit $\lambda$ as $R\rightarrow +\infty$. Therefore, this explicit example tells us that,
if $\gamma$ is near to 2, $c$ is sufficiently large,
and $c^2\Lambda/(4\pi GA_1)$ is near to $\lambda >1$, then $P(r)$ of the prolongation of the germ, which
satisfies $dP/dr>0$ as $0<r \ll 1$, become decreasing. \\

In the definition of `monotone-short' solutions of the
Tolman-Oppenheimer-Volkoff-de Sitter equation we have
required that $\kappa_+>0$ and $Q_+>0$. Let us spend few words concerning these conditions.

Consider a solution $(m,P)=(m(r), P(r)), 0<r<r_+,$ in $\mathcal{D}$ such that $dP/dr<0$,
which requires $Q(r, m(r), P(r))>0$
 for $0<r<r_+$, and suppose $P(r)\rightarrow 0$ as
$r\rightarrow r_+-0$, with $r_+$ being finite.

When we are concerned with the Tolman-Oppenheimer-Volkoff equation (\ref{TOV8})
with $\Lambda=0$, the condition $\kappa_+>0$ follows automatically. Proof can be found in \cite{2015KJM}. Of course, if $\Lambda=0$,  then $Q_+=Gm_+>0$ a priori.

However if $\Lambda >0$ it seems that we cannot exclude the possibility
that $\kappa_+=0$ a priori. Generally speaking, since $\kappa>0$ in $\mathcal{D}$,
we have $\kappa_+\geq 0$. Suppose $\kappa_+=0$. Then
$$\frac{2Gm_+}{c^2r_+}=1-\frac{\Lambda}{3}r_+^2,
$$
and, since $\kappa>0$ for $r<r_+$ and $\kappa_+=0$, we see
\begin{align*}
\kappa_+'&:=\frac{d\kappa}{dr}\Big|_{r=r_+-0}=\lim_{r\rightarrow r_+-0}-\frac{2G}{c^2}4\pi \rho+
\frac{2Gm}{c^2}\frac{1}{r^2}-\frac{2}{3}\Lambda r \\
&=\frac{2Gm_+}{c^2}\frac{1}{r_+^2}
-\frac{2}{3}\Lambda r_+  =\frac{1}{r_+}(1-\Lambda r_+^2) \leq 0.
\end{align*}
On the other hand, we have
\begin{align*}
Q_+&=Gm_+-\frac{c^2\Lambda}{3}r_+^3 =\frac{c^2r_+}{2}(1-\Lambda r_+^2)\\
& \geq 0,
\end{align*}
since $Q>0$ for $r<r_+$. Thus it should be the case that
$1-\Lambda r_+^2=0$ and $Q_+=0$. In other words, $\kappa_+=0$
requires $Q_+=0$ and $\Lambda r_+^2=1$. This is very non-generic situation probably hard to occur,
but at the moment we have no reason to exclude this possibility.

\section{Metric on the vacuum region}

Suppose that we have fixed a solution $(m(r),P(r)),
0<r<r_+,$ of the Tolman-Oppenheimer-Volkoff-de Sitter equation (\ref{TOV1})
which is monotone-short. Then we have the metric
$$ds^2=\kappa_+e^{-2u/c^2}c^2dt^2-
\frac{1}{\kappa}dr^2-r^2d\omega^2 
$$
on $0\leq r <r_+$, where
$$d\omega^2=d\theta^2+\sin^2\theta d\phi^2.$$
We should continue this metric to the exterior
vacuum domain $r\geq r_+$. Naturally, keeping in mind the Birkhoff
theorem, we should take the Schwarzschild-de Sitter
metric
$$ds^2=\Big(1-\frac{2Gm_+}{c^2r}-\frac{\Lambda}{3}r^2\Big)c^2dt^2
-\Big(1-\frac{2Gm_+}{c^2r}-\frac{\Lambda}{3}r^2\Big)^{-1}dr^2
-r^2d\omega^2
$$
on $r\geq r_+$. As the whole we take
$$ds^2=g_{00}c^2dt^2-g_{11}dr^2-r^2d\omega^2,
$$
where
\begin{align*}
g_{00}&=\begin{cases}
\kappa_+e^{-2u(r)/c^2} & \quad (0\leq r <r_+) \\
\displaystyle 1-\frac{2Gm_+}{c^2r}-\frac{\Lambda}{3}r^2 &\quad (r_+\leq r <r_E)
\end{cases},\\
-g_{11}&=\Big(1-\frac{2G\tilde{m}(r)}{c^2r}
-\frac{\Lambda}{3}r^2\Big)^{-1} \quad (0\leq r < r_E),
\end{align*}
with
$$\tilde{m}(r)=\begin{cases}
m(r) & \quad (0\leq r <r_+) \\
m_+ &\quad (r_+\leq r<r_E).
\end{cases}
$$
Here the constants $r_E, r_I, (0<r_I<r_E<+\infty),$ are the values of $r$ of the so called `cosmological horizon', `black hole horizon', that is, 
$\kappa(r,m_+)>0$
if and only if $r_I<r<r_E$. 
In other words, we have
$$\kappa(r,m_+)=\frac{\Lambda}{3r}(r-r_I)(r_E-r)(r+r_I+r_E).
$$
See \cite{Brady}. But this situation is possible only
if
\begin{equation}
\sqrt{\Lambda} <\frac{c^2}{3Gm_+}. \label{B}
\end{equation}
If (\ref{B}) does not hold, then $\kappa(r,m_+)
\leq 0$ for all $r>0$.
However, since we are
supposing  $\kappa_+=\kappa(r_+,m_+)>0$, the condition (\ref{B})
is supposed to hold and we have 
$r_I<r_+<r_E$.
Let us discuss the regularity of this patched metric.\\

First we observe the regularity of $u(r)$. We claim

\begin{Proposition}
The function $u(r)$ is of class $C^2$ in a neighborhood of
$r_+$ and
$$u(r)=B(r_+-r)(1+O(r_+-r))
$$
as $r\rightarrow r_+-0$, with $B:=Q_+/r_+^2\kappa_+$.
Hence
$\rho(r)$ is of class $C^1$ and
$$\rho(r)=\Big(\frac{(\gamma-1)B}{\gamma A}\Big)^{\frac{1}{\gamma-1}}
(r_+-r)^{\frac{1}{\gamma-1}}(1+O(r_+-r)).
$$
\end{Proposition}

Proof.  Since $u(r)$ satisfies the equation (\ref{TOV5}), whose right-hand 
side is a $C^1$-function of $(r,m,u)$ near $(r_+, m_+, 0)$
thanks to $\displaystyle \mu=\frac{1}{\gamma-1}>1$. Recall that $\kappa_+>0$.
Therefore the continuous solution $u(r)$ turns out to be
of class $C^2$ and
$$\frac{du}{dr}\Big|_{r=r_+-0}=-\frac{Q+}{r_+^2\kappa_+}=-B.
$$
This completes the proof.\\

Now we are going to see the regularity of $g_{00}, g_{11}$.
Since
$$
\frac{d}{dr}\tilde{m}(r)=
\begin{cases}
4\pi r^2\rho(r) &\quad (r<r_+) \\
0 &\quad (r_+\leq r <r_E)
\end{cases}
$$
is of class $C^1$, $\tilde{m}(r)$ is of class $C^2$.
Therefore $g_{11}$ is twice continuously differentiable
across $r=r_+$.

Since $u$ vanishes at $r=r_+-0$, $g_{00}$ is continuous thanks to the
definition of $\kappa_+$. We see
$$\frac{d}{dr}g_{00}\Big|_{r=r_+-0}=-\frac{2\kappa_+}{c^2}\frac{du}{dr}
\Big|_{r=r_+-0}=\frac{2Q_+}{c^2r_+^2}
$$
and
$$\frac{d}{dr}g_{00}\Big|_{r=r_++0}=
\Big(\frac{2Gm_+}{c^2r^2}-
\frac{2\Lambda}{3}r\Big)_{r=r_+}=
\frac{2Q_+}{c^2r_+^2}.
$$
Therefore $g_{00}$ is continuously differentiable. We have
$$\frac{d^2}{dr^2}g_{00}\Big|_{r=r_+-0}=
\frac{4\kappa_+}{c^4}\Big(\frac{du}{dr}\Big)^2_{r=r_+-0}
-\frac{2\kappa_+}{c^2}\Big(\frac{d^2u}{dr^2}\Big)_{r=r_+-0}.
$$
But by a tedious calculation we have
$$\frac{d^2u}{dr^2}\Big|_{r=r_+-0}=
\frac{c^2\Lambda}{\kappa_+}+\frac{2Q_+}{r_+^3\kappa_+}+
\frac{2(Q_+)^2}{c^2r_+^4\kappa_+^2}.
$$
This can be derived by differentiating the right-hand 
side of the equation for $du/dr$,
that is, the second equation of (\ref{TOV5}). Therefore we see
$$\frac{d^2}{dr^2}g_{00}\Big|_{r=r_+-0}=
\frac{d^2}{dr^2}g_{00}\Big|_{r=r_++0}=
-\frac{4Q_+}{c^2r_+^3}-2\Lambda.
$$
Hence $g_{00}$ is twice continuously differentiable
across $r=r_+$. Summing up, we have
\begin{Theorem}
Given a monotone-short solution of the 
Tolman-Oppenheimer-Volkoff-de Sitter equation (\ref{TOV1}), we can extend the interior metric
to the exterior Schwarzschild-de Sitter metric on the vacuum
region with twice continuous differentiability.
\end{Theorem}

\section{Analytical property of the vacuum boundary}

Let us observe the analytical property of
a monotone-short solution $(m(r), P(r))$, $0<r<r_+,$
of the Tolman-Oppenheimer-Volkoff-de Sitter equation
(\ref{TOV1}). Proposition 2 tells us that
the associated $u(r)$ belongs to $C^2([0,r_+])$ and
$$u(r)=B(r_+-r)(1+O(r_+-r))$$
as $r\rightarrow r_+-0$, where $B=Q_+/r_+^2\kappa_+$. 
Moreover we claim

\begin{Theorem}
Any monotone-short solution $u(r), 0<r<r_+,$ of (\ref{TOV1})
enjoys the behavior at $r=r_+-0$ such that
$$
u(r)=B(r_+-r)(1+[r_+-r, (r_+-r)^{\frac{\gamma}{\gamma-1}}]_{1}), $$
therefore
$$
\rho(r)=
\Big(\frac{(\gamma-1)B}{\gamma A}\Big)^{\frac{1}{\gamma-1}}
(r_+-r)^{\frac{1}{\gamma-1}}
(1+[r_+-r, (r_+-r)^{\frac{\gamma}{\gamma-1}}]_{1}). $$
Here $[X_1, X_2]_{1}$ stands for a convergent double
power series of the form 
$$\sum_{k_1+ k_2\geq 1}a_{k_1k_2}X_1^{k_1}X_2^{k_2}.$$
\end{Theorem}

Proof.
Let us denote $\displaystyle \mu:=\frac{1}{\gamma-1}$ so that
$\displaystyle \frac{\gamma}{\gamma-1}=\mu +1$.

First suppose that $\mu$ is an integer. Then proof  is easy. 
In fact $(m(r), u(r))$ satisfies at least on $0<r<r_+$ the system of equations
\begin{align}
\frac{dm}{dr}&=
4\pi r^2A_1u^{\mu}\Omega_{\rho}(u/c^2), \label{TOV01} \\
\frac{du}{dr}&=
-\frac{G\Big(\displaystyle m+
\frac{4\pi}{c^2}r^3AA_1^{\gamma}
u^{\mu+1}\Omega_P(u/c^2)\Big)-\frac{c^2\Lambda}{3}r^3}{r^2\Big(\displaystyle 1-\frac{2Gm}{c^2r}-\frac{\Lambda}{3}r^2\Big)}. \nonumber
\end{align}
and $(m(r), u(r)) \rightarrow (m_+, 0)$ as $r\rightarrow r_+-0$.
But, since $\mu$ is supposed to be an integer, the right-hand side of the system (\ref{TOV01}) is analytic
function of $(r, m, u)$ in a neighborhood of
$(r_+, m_+, 0)$. This guarantees that
$m(r), u(r)$ admit analytic prolongations beyond $r=r_+$ to the right,
and completes the proof. Of course this analytic prolongation is different from
the $C^2$-prolongation as a solution of (\ref{TOV5}), since 
$u^{\mu} \not= (u_{\sharp})^{\mu}(=0)$ for $u <0$.

Now suppose that $\mu$ is not an integer. Since
$u(r), 0<r<r_+, $ is monotone decreasing, it has the inverse function
$r=r(u)$ defined on $0<u<u_c$ such that $r(u)\rightarrow r_+$
as $u\rightarrow +0$. Then we have a solution
$(m,r)=(m(u), r(u))$ of the system of equations
\begin{subequations}
\begin{eqnarray}
\frac{dm}{du}&=&-4\pi r^4\Big(1-\frac{2Gm}{c^2r}-
\frac{\Lambda}{3}r^2\Big)\cdot Q^{-1}\cdot A_1u^{\mu}\Omega_{\rho}(u/c^2),  \label{TOV02a}\\
\frac{dr}{du}&=&-r^2\Big(1-\frac{2Gm}{c^2r}-\frac{\Lambda}{3}r^2\Big)\cdot Q^{-1}, \label{TOV02b}
\end{eqnarray}
\end{subequations}
where $$Q=G\Big(m+\frac{4\pi}{c^2}r^3AA_1^{\gamma}
u^{\mu+1}\Omega_P(u/c^2)\Big)-\frac{c^2\Lambda}{3}r^3.$$
Since $(m(u), r(u)) \rightarrow (m_+, r_+)$ and
$Q\rightarrow Q_+>0$ as $u\rightarrow +0$ and
the right-hand sides of (\ref{TOV02a})(\ref{TOV02b}) are analytic functions of
$u, u^{\mu}, m, r$ on a neighborhood of
$(0,0, m_+, r_+)$, we can apply the following Lemma
in order to get
\begin{align*}
m(u)&=m_++u[u, u^{\mu}]_0, \\
r(u)&=r_++u[u, u^{\mu}]_0.
\end{align*}
Here $[\cdot, \cdot]_0$ denotes a convergent double power series.

\begin{Lemma}
Let $\mu >1$ and 
$f^{\alpha}(x, x^{\mu}, y_1, y_2), \alpha=1,2, $
be analytic functions of
$x, x^{\mu}, y_1, y_2$ on a neighborhood of $(0,0,0,0)$.
Let $(y_1(x), y_2(x)), 0\leq x\leq \delta,$ be the solution of the problem
\begin{equation}\frac{dy_{\alpha}}{dx}=f^{\alpha}(x,x^{\mu},y_1,y_2),\quad
y_{\alpha}|_{x=0}=0, \qquad \alpha=1,2. \label{EqLemma}
\end{equation}
Then there are analytic functions $\varphi^{\alpha}$ of $x, x^{\mu}$
on a neighborhood of $(0,0)$ such that
$y_{\alpha}(x)=x\varphi^{\alpha}(x,x^{\mu})$ for $0<x\ll 1$.
\end{Lemma}

A proof of this Lemma will be sketched in the Appendix.\\

Since $dm/du \sim -C u^{\mu}$, with $C=
4\pi r_+^4(\kappa_+/Q_+)A_1$, and $du/dr\rightarrow -B$, we have
\begin{align*}
m&=m_+-Cu^{\mu+1}+\sum_{n\geq 2}m_{1n}u^{\mu n+1}+
\sum_{n\geq 0, l\geq 2}m_{ln}u^{\mu n+l}, \\
r&=r_+-\frac{1}{B}u+\sum_{n\geq 1}c_{1n}u^{\mu n+1}
+\sum_{n\geq 0, l\geq 2}c_{ln}u^{\mu n+l}.
\end{align*}
If $m_{1n}\not=0$ for $\exists n\geq 2$, then $dm/du$ would
contain the term $u^{\mu n}$ with $n\geq 2$. However it is impossible,
since the right-hand side of (\ref{TOV02a}) cannot contain such a term. Therefore $m_{1n}=0$ for $\forall n\geq 2$.
If $c_{1n}\not=0$ for $\exists n\geq 1$, then $dr/du$ would contain the term
$u^{\mu n}$. However it is impossible, since
the right-hand side of the equation (\ref{TOV02b}) cannot contain such a term. Therefore $c_{1n}=0$ for $\forall n\geq 1$. Thus we have
\begin{align*}
m&=m_+-Cu^{\mu+1}+
\sum_{n\geq 0, l\geq 2}m_{ln}u^{\mu n+l}, \\
r&=r_+-\frac{1}{B}u+\sum_{n\geq 0, l\geq 2}c_{ln}u^{\mu n+l}
\end{align*}
Moreover we can show that $m_{ln}=c_{ln}=0$ for $2\leq l\leq n$
by induction on $l$. In fact, fix $n\geq 2$. Then $m_{2n}=c_{2n}=0$,
since, otherwise, $dm/du$, therefore the right-hand side of (\ref{TOV02a}), or
$dr/du$, therefore the right-hand side of 
(\ref{TOV02b}), would contain the term $u^{\mu n+1}$ with $n\geq 2$,
which is impossible. Therefore $m_{2n}=c_{2n}=0$. Consider $3\leq l\leq n$.
Assume $m_{n',l'}=c_{n',l'}=0$ for $2\leq l'\leq n'$ with $n'\leq n, l'\leq l-1$. If $m_{nl}\not=0$, or, $c_{nl}\not=0$, then
$dm/du$, or $dr/du$ would contain the term $u^{\mu n+l-1}$,
which is impossible by the induction assumption. Therefore
$m_{nl}=c_{nl}=0$ for $2\leq l\leq n$. This implies that
\begin{align*}
r&=r_+-\frac{1}{B}u+\sum_{n\geq 0, l\geq n+1, l\geq 2}c_{ln}
u^{(\mu+1)n+l-n} \\
&=r_+-\frac{1}{B}u(1+[u, u^{\mu+1}]_1).
\end{align*}
The inverse function $u=u(r)$ then clearly enjoys an expansion of the form
$$u=B(r_+-r)(1+[r_+-r, (r_+-r)^{\mu+1}]_1).$$
This completes the proof of Theorem 4.\\

{\bf\Large Acknowledgment} 

The author would like to express his sincere thanks to
Professor Cheng-Hsiung Hsu (National Central University, Taiwan) who kindly read the original manuscript carefully and indicated typos.
\\

{\bf\Large Appendix}\\

Le us sketch a proof of Lemma 1.

First we assume that $\mu$ is not an integer but a rational
number, say, $\mu=q/p, p, q \in \mathbb{N}, p\geq 2$ and
$p,q$ are relatively prime. Note that a function given by a convergent power series 
$$\varphi(x)=\sum \tilde{c}_{ij}x^i(x^{\mu})^j,\qquad |\tilde{c}_{ij}|\leq 
\frac{\tilde{M}}{\tilde{\delta}^{i+j}}, \quad (\tilde{\delta}<1), $$
can be rewritten as
$$\varphi(x)=\sum_{l}\sum_{n=0}^{p-1} c_{l n}x^{\mu n+l},\qquad |c_{ln}|\leq \frac{M}{\delta^{l}} \quad \mbox{for}\quad 0\leq n\leq p-1,$$
where
$$c_{ln}=\sum\{\tilde{c}_{ij}\  |\  i+qJ=l, j=pJ+n, \exists J\in\mathbb{N}\},$$
and
$M=\tilde{M}/(e \tilde{\delta}^p), \delta=\tilde{\delta}/e$.
This rewriting is necessary, since not $\tilde{c}_{ij}$'s but
$c_{ln}$'s can be uniquely determined for the given function $\varphi(x)$.

In fact first we note that $\mu n$ cannot be an integer for $n=1,\cdots, p-1$. (Proof:
Let us deduce a contradiction supposing that $nq/p$ is an integer.
We can assume $q<p$, by, if necessary, replacing $q$ by $q':=q-[q/p]p$. Since $nq/p <n$, we see that
$nq/p$ is either $ 1,\cdots,$ or $n-1$, therefore
$q/p$ is either $1/n,\cdots,$ or $(n-1)/n$.
Hence $p$ is a divisor of $n$, a fortiori, $p\leq n$, a contradiction to
$n\leq p-1$, QED.)
Hence $\mu n+l=\mu n'+l', n,n',l,l'\in\mathbb{N}, 0\leq n, n'\leq p-1,$ implies $n=n', l=l'$.
Then we have a unique numbering $(n_k, l_k)_{k\in\mathbb{N}}$ of $(n,l)$'s such that 
$\mu n_k+l_k < \mu n_{k+1}+l_{k+1}$. By induction on $k$ we can deduce $c_{l_kn_k}=0$
for $\forall k$ from $\sum c_{ln}x^{\mu n+l}=0\   \forall x$. This means the uniqueness of the coefficients $c_{ln}$ in the above expansion of $\varphi(x)$. 

Anyway suppose
$$f^{\alpha}(x, x^{\mu}, y_1, y_2)=\sum\sum_{n=0}^{p-1}
a_{lnk_1k_2}^{\alpha}x^{\mu n+l}y_1^{k_1}y_2^{k_2},$$
with
$$\Big|a_{lnk_1k_2}^{\alpha}\Big|\leq \frac{M}{\delta^{l+k_1+k_2}}
\qquad (0\leq n\leq p-1) $$
and put
\begin{align*}
F(x,y_1,y_2)&=\sum\frac{M}{\delta^{l+k_1+k_2}}
\sum_{n=0}^{p-1}x^{\mu n+l}y_1^{k_1}y_2^{k_2} \\
&=\frac{M}{1-x/\delta}
\frac{1-x^{\mu p}}{1-x^{\mu}}
\frac{1}{1-y_1/\delta}\frac{1}{1-y_2/\delta}.
\end{align*}
Then the problem 
$$\frac{dY}{dx}=F(x,Y, Y),\qquad Y|_{x=0}=0$$
has a solution of the form
\begin{align*}
Y&=Mx(1+[x, x^{\mu}]_1) \\
&=\sum_l\sum_{n=0}^{p-1}C_{ln}x^{\mu n+l}, \qquad
0\leq C_{ln}\leq \frac{M'}{(\delta')^l}.
\end{align*}
On the other hand (\ref{EqLemma}) has a formal power series solution
$$y_{\alpha}=\sum_l\sum_{n=0}^{p-1}c_{ln}^{\alpha}x^{\mu n+l}, $$
where the coefficients $c_{ln}^{\alpha}$'s are determined by a recursive
formula
\begin{align*} 
c_{0n}^{\alpha}&=0,\qquad c_{l+1, n}^{\alpha}=\frac{1}{l+1+\mu n}b_{ln}^{\alpha}, \\
b_{LR}^{\alpha}&=\sum a_{lnk_1k_2}^{\alpha}c_{l'(1)n'(1)}^1\cdots
c_{l'(k_1)n'(k_1)}^1  c_{l''(1)n''(1)}^2\cdots c_{l''(k_2)n''(k_2)}^2.
\end{align*}
Here the summation in the definition of $b_{LR}^{\alpha}$ is taken
over $$L=qJ+l+l'(1)+\cdots +l'(k_1)+l''(1)+\cdots
+l''(k_2)
$$
with $l'(1),\cdots, l''(1), \cdots \geq 1$ and $J\in\mathbb{N}$,
and
$$pJ+R=n+n'(1)+\cdots +n'(k_1)+
n''(1)+\cdots+n''(k_2).$$
Then it can be shown inductively that
$|c_{ln}^{\alpha}|\leq C_{ln}$, which implies the convergence of the formal power series solution. This completes the proof.

A proof by a similar and easier majorant argument can be done when $\mu$ is an irrational
number. Let us omit the repetition.

\vspace{15mm}
Professor Emeritus at Yamaguchi University

e-mail: makino@yamaguchi-u.ac.jp

\end{document}